\documentclass[11pt]{amsart}

\usepackage{amsmath,amsthm,amsfonts,amssymb,bm, graphicx}
\usepackage{xcolor}
 \usepackage{soul}
   \usepackage[parfill]{parskip} 
\usepackage{hyperref, mathrsfs}
 \numberwithin{equation}{section}
   \theoremstyle{plain}
 
\newtheorem*{con* A}{Conjecture A}
\newtheorem*{con* B}{Conjecture B}

\newtheorem{thm*}{Theorem}
\newtheorem*{con*}{Conjecture}
\newtheorem*{lem*}{Lemma}

\theoremstyle{definition}

\theoremstyle{definition}
\newtheorem*{Acknowledgements}{Acknowledgement}

\renewcommand{\geq}{\geqslant}
\renewcommand{\leq}{\leqslant}

\newcommand{\eps}{\varepsilon}

\newcommand\be{\begin{equation}}
\newcommand\ee{\end{equation}}
\newcommand\bes{\begin{equation*}}
\newcommand\ees{\end{equation*}}
 
 \newcommand{\est}[1]{\begin{equation*}\begin{split}#1\end{split}\end{equation*}}

\newcommand{\R}{\mathbb{R}}

\begin{document}


\title{The $\theta=\infty$ conjecture implies the Riemann hypothesis}
\author {Sandro Bettin} 
\author{Steven M. Gonek}

\address{DIMA - Dipartimento di Matematica, Via Dodecaneso, 35 \\
16146 Genova - ITALY} 
\email{bettin@dima.unige.it}

\address{Department of Mathematics, University of Rochester, Rochester, NY 14627}
\email{gonek@math.rochester.edu}

\thanks{Work of the first author was partially supported by NSF grant DMS-1200582.}



\begin{abstract}
We show that the $\theta=\infty$ conjecture implies the Riemann hypothesis.
\end{abstract}

\maketitle

\section{Introduction}
Since the work of Levinson~\cite{Lev}, it has been known that one can obtain lower bounds for the proportion of zeros of the Riemann zeta-function on the critical line by computing upper bounds for the mollified second moment
\be\label{def I}
I_N(T_1,T_2) :=\int_{T_1}^{T_2} |M_{N}(\tfrac12+it)|^2|\zeta(\tfrac12+it)|^2\,dt,
\ee
where $M_N(s)$ is a mollifier   roughly of the form
$$
M_N(s):=\sum_{n\leq N}\frac{\mu(n)}{n^{s}}\Big(1-\frac{\log n}{\log N}\Big)  
$$
with $N\geq 2$   an  integer.
Levinson~\cite{Lev} computed the asymptotic  formula
\be\label{ff}
\lim_{T\rightarrow\infty}\frac{ I_{T^{\theta} }(0,T)}T=1+\frac1{\theta}
 \ee
for $0<\theta<\frac12$, and used this result to deduce that $\kappa>\frac13$, where
$$
\kappa:= \frac{\#\{\rho \mid\zeta(\rho)=0,\ 0<\Im\,\rho<T,\ \Re\,\rho=\frac12\}}{\#\{\rho \mid\zeta(\rho)=0,\ 0<\Im\,\rho<T\} }
$$ 
is the proportion of the non-trivial zeros of $\zeta(s)$ that lie on the critical line.
Conrey~\cite{Con} later proved that~\eqref{ff} (with a slightly different mollifier)   remains valid   for $\theta<\frac47$, and thereby deduced that $\kappa>\frac25$. 

Initially it was believed (see \cite{Far}) that~\eqref{ff} does not hold when 
$\theta>1$. However,   Farmer~\cite{Far} produced a  heuristic argument  
suggesting  that it   holds for every $\theta>0$, and  called this the  
 ``$\theta=\infty$ conjecture''. Moreover, he proved that this conjecture implies that
  $\kappa=1$, in other words, that $100\%$ of the non-trivial zeros of $\zeta(s)$ lie on the critical line.  He also argued that a slightly stronger form of the conjecture implies Montgomery's pair correlation conjecture.
More recently, Radziwi\l\l~\cite{Rad} showed that, as $\theta\rightarrow\infty$, $M_{T^{\theta}}(t)$ is essentially the best possible mollifier of length $T^{\theta}$ for $\zeta(s)$. In particular, his work implies that Levinson's method can give $\kappa=1$ only if it is used with mollifiers of length $T^\theta$, where $\theta$ is arbitrarily large. 

The purpose of this note is to show that the $\theta=\infty$ conjecture actually implies the Riemann hypothesis. Indeed, we show that
even an upper bound of the form $I_N(0, T)\ll T^{1+\eps}$   for some 
$\theta>1$ and all $N$ in the range   $2\leq N\leq T^{\theta}$ implies a zero-free region for the zeta-function of the form $\Re\, s >1-\delta$ for some $\delta>0$ depending on $\theta$; in other words,    a quasi-Riemann hypothesis.

\begin{thm*}\label{t1}
Let $\theta>0$ and assume   that for every $\eps>0$ we have $I_N(0,T )\ll_\eps T^{1+\eps}$ for  $N$ in the range $2\leq N\leq T^{\theta}$. Then $\zeta(s)$ has no zeros in the half-plane $\Re\,s> \frac12+ \frac{1}{2\theta}$. In particular, if $I_N(0,T)\ll_\eps T^{1+\eps}$ for   $2\leq N\leq T^{\theta}$ with $\theta$ arbitrarily large, then the Riemann hypothesis is true. 
\end{thm*}

In a number of  recent works on mean values of $L$-functions in the $t$-aspect, the integral  is taken over   $[T,2T]$ rather than over $[0,T]$. 
Thus,  it is natural to ask whether one can   obtain a version of Theorem~\ref{t1} for the interval $[T,2T]$. Usually there is no difficulty in passing from one    interval to the other.
In our case,  however,  the problem for $[T,2T]$ is more subtle because one needs an $\Omega$-result for $M_N(t)$ that is uniform in $t$.
Using ideas from~\cite{Pin} and~\cite{GGL}, we prove the following. 

\begin{thm*}\label{t2}
Let $\theta>0$ and assume that for every $\eps>0$ we have $I_N(T,2T)\ll_\eps T^{1+\eps}$ for  $N$ in the range $2\leq N\leq T^{\theta}$. Then $\zeta(s)$ has no zeros in the half-plane $\Re\,s > \frac12+\frac{2}{\theta}$. In particular, if $I_N(T,2T)\ll_\eps T^{1+\eps}$ for   $2\leq N\leq T^{\theta}$ with $\theta$ arbitrarily large, then the Riemann hypothesis is true.
\end{thm*}

Notice that  Theorem~\ref{t2} only implies a  quasi-Riemann hypothesis when 
$\theta>4$, so in this respect it is weaker than Theorem~\ref{t1}. However, Theorem~\ref{t2}, whose proof is more difficult than that of Theorem~\ref{t1}, is in a certain sense best possible. If, for example, one assumes 
that $\zeta(s)$ has a unique simple zero $\rho_0=\beta_0+i\gamma_0$ such that $\gamma_0>0$ and $\beta_0>\frac12$,   
one can show    that
\est{
I_N(T,2T)=c_1\frac{N^{2\beta_0-1}}{T^{3}}\frac{\log T}{\log^2N}\Big(1+\Re\Big(N^{2i\gamma_0}\frac{|\zeta'(\rho_0)|^2}{\zeta'(\rho_0)^2}\Big)+o(1)\Big) +O\Big(T^{1+\eps}+\frac{N^{\beta_0-\frac12+\eps}}{T}\Big)
}
for some   constant $c_1>0$, as $T\rightarrow\infty$, and this is consistent 
with the assumption $I_{T^{\theta}}(T,2T)\ll T^{1+\eps}$ if $\theta<4$. For the sake of comparison, we note that with the same   zero configuration one has
\est{
I_N(0,T)=\frac{N^{2\beta_0-1}}{\log^2 N}(C(N)+o(1))+O(T^{1+\eps}+N^{\beta_0-\frac12+\eps}T^\eps)
}
for some positive function $C(N)$ bounded away from $0$, so that $I_{T^{\theta}}(0,T)\ll T^{1+\eps}$ implies $\beta_0\leq\frac12+\frac1{2\theta}$, 
which is consistent with Theorem~\ref{t1}.
\begin{Acknowledgements}
The first author would like to thank Brian Conrey and Jon Keating for bringing this problem to his attention. 
Both authors wish to thank the organizers of the workshop ``Computational Aspects of L-functions'' and ICERM for providing an excellent environment for collaboration.
\end{Acknowledgements}

\section{Proof of the Theorems}

We will prove Theorem~\ref{t1} and~Theorem~\ref{t2} at the same time. It should be pointed out, however, that an easier argument would suffice for the former.

We begin by extending our earlier definition of $M_N(s)$ slightly by writing
\be
M_x(s)\log x 
=\sum_{n\leq x}\frac{\mu(n)}{n^{s}} \log (x/n)  
\ee
for $x>0$ (with $M_1(s):=0$). Notice that the right-hand side is zero when $0<x\leq 1$ and that 
this also allows us to extend  the definition of $I_N(T_1, T_2)$ in \eqref{def I}
to $I_x(T_1,T_2)$. 
Now, for   $t\in\R$ we   have
 \be\notag
 M_x(\tfrac12+it)\log x=\frac1{2\pi i  }\int_{1-i\infty}^{1+i\infty}\frac{x^z}{\zeta(\frac12+it+z)}\,\frac{dz}{z^2}.
 \ee
Thus, by Mellin inversion we see that
\est{
H_t(w)&:=\int_{1}^{\infty}M_x(\tfrac12+it)(\log x) x^{-w}\,dx=\frac1{(w-1)^2\zeta(w-\frac12+it)}\\
}
for $\Re\, w>\frac32$.
Next, assuming that $\rho_0= \beta_0+i\gamma_0$ is a fixed zero of $\zeta(w)$ with 
$\beta_0 \geq 1/2$, we define
\begin{equation}\notag
\begin{split} 
G_t(w) :=\frac{(w-1)^2(w-\frac32+it)\zeta(w-\frac12+it)}{(w+1)^2(w-\frac12+it-\rho_0)  (w+it+1)^4}.
\end{split}
\end{equation}
In the half-plane $\Re\, w \geq 0$,   $G_t(w)$ is holomorphic and satisfies $G_t(w)\ll(1+ |w+it|)^{-\frac52}$. Thus,   setting
\begin{equation}\notag
\begin{split} 
g_t(u) = \frac{1}{2\pi i} \int_{3-i\infty}^{3+i\infty} G_t(w)u^{-w} dw 
\end{split}
\end{equation}
for   $u>0$, we have 
\be\label{gbound}
g_t(u)=\begin{cases}0 & \text{if $u>1$},\\
O(1) & \text{if $0\leq u\leq 1$},
       \end{cases}
\ee
as can be seen by moving the line of integration to $\Re\, w =+\infty$ when $u>1$, 
and to $\Re\, w=0$ when $0\leq u\leq1$. 

Now consider the integral
\be\label{def J}
J_t(x):=\frac{1}{2\pi i} \int_{3-i\infty}^{3+i\infty} G_t(w)H_t(w)x^{w} dw=\frac{1}{2\pi i} \int_{3-i\infty}^{3+i\infty} \frac{(w-\frac32+it)x^{w}}{(w+1)^2(w-\frac12+it-\rho_0)  (w+it+1)^4}  dw, 
\ee 
where, from this point on, we assume that $x\geq 2$.
On the one hand,   by the convolution formula for   products of Mellin transforms,  and since $M_{y}(\tfrac12+it)\log y =0$ when $0< y \leq 1$,
\bes
J_t(x)=\int_1^{\infty} M_{y}(\tfrac12+it)(\log y)g_t (y/x)\, dy.
\ees
Thus, by ~\eqref{gbound}, 
\be\label{1I}
J_t(x)\ll \int_1^{x} |M_{y}(\tfrac12+it)| \log y\, dy
\ee
for $x\geq 2$.
On the other hand, moving the line of integration in \eqref{def J} to $\Re\, w=0$, we see that
\be\label{2I}
J_t(x)=\frac{1}{2\pi i} \int_{-i\infty}^{+i\infty} G_t(w)H_t(w)x^{w} dw+\frac{x^{\rho_0+\frac12-it}(\rho_0-1)}{(\frac32+\rho_0-it)^2(\rho_0+\frac32)^4}.
\ee
The integral on the right  is $O(1)$ since $H_t(w)G_t(w)\ll(1+ |w|)^{-2}$ for $\Re\, w=0$. Thus, from~\eqref{1I} and~\eqref{2I} we deduce that
\bes
\frac{x^{\beta_0+\frac12}}{(1+|t|)^2} +1  \ll   \int_{1}^{x} |M_{y}(\tfrac12+it)| \log y \,dy.
\ees
It   follows from the Cauchy-Schwarz inequality that 
\bes
\frac{x^{2\beta_0}}{(1+|t|)^4}+ \frac{1}{x }  \ll \int_1^{x} |M_{y}(\tfrac12+it)|^2  \log^2 y\, dy
\ees
for $x\geq 2$.
Multiplying both sides by $|\zeta(\frac12+it)|^2$ and integrating with respect to $t$ over the interval 
$[T_1,T_2]$, where $0\leq T_1\leq T_2/2$, we obtain
\bes
\begin{split}
 \int_{T_1}^{T_2}|\zeta(\tfrac12+it)|^2
\bigg(\frac{x^{2\beta_0}}{(1+t)^4} +\frac{1}{x}\bigg) \,dt
 \ll  &\int_1^{x} \log^2 y \int_{T_1}^{T_2} |M_{y}(\tfrac12+it)\zeta(\tfrac12+it)|^2 \,dt dy\\
\leq& \log^2 x\int_1^{x} I_y(T_1,T_2)\,  dy.
 \end{split}
\ees
Now $\int_{T_1}^{T_2}|\zeta(\frac12+it)|^2\,dt \gg T_2\log( T_2+2)$ for $0\leq T_1\leq T_2/2$,  so 
\bes
\frac{x^{ 2\beta_0 }\log(T_1+2)}{|1+T_1|^3} 
+\frac{T_2\log (T_2+2)}{x}
\ll \log^2x\int_1^{x} I_y(T_1,T_2)\,  dy.
\ees
Thus, if $I_N(0,T)\ll_\eps T^{1+\eps}$ holds for   $2\leq N\leq T^{\theta}$ and for every $\eps>0$, then taking $T_1=0$, $T_2=T$,  and $x=T^{\theta}$, we obtain
\bes
T^{ 2\beta_0\theta}\ll_\eps T^{1+\eps+\theta}.
\ees
Letting $T\rightarrow\infty$ and letting $\eps>0$ be sufficiently small, we obtain $\beta_0\leq \frac12+\frac1{2\theta}$, as claimed in Theorem~\ref{t1}. Theorem~\ref{t2} follows in the same way on taking $T_1=T$ and $T_2=2T$.

\appendix




\begin{thebibliography}{0}

\bibitem{Con}
Conrey, J.B. \emph{More than two fifths of the zeros of the Riemann zeta function are on the critical line}. J. Reine Angew. Math. 399 (1989), 1--26.
 
\bibitem{Far}
Farmer, D.W. \emph{Long mollifiers of the Riemann zeta-function}. Mathematika 40 (1993), no. 1, 71--87. 

\bibitem{GGL} 
Gonek, S. M., Graham, S. W., and Lee Y. 
\emph{A Generalized Lindel\"of Hypothesis}, unpublished manuscript.

\bibitem{Lev}
Levinson, N. \emph{More than one third of zeros of Riemann's zeta-function are on $\sigma=1/2$}. Advances in Math. 13 (1974), 383--436. 

\bibitem{Pin}
Pintz, J. Oscillatory properties of $M(x)=\sum _{n\leq x}\mu(n)$, I. Acta Arith. 42 (1982/83), no. 1, 49--55.

\bibitem{Rad}
Radziwi\l\l, M. \emph{Limitations to mollifying $\zeta(s)$}, preprint; arxiv math.NT/1207.6583.


\end{thebibliography}
\end{document}